\documentclass[a4paper,10pt,twoside]{amsart}

\usepackage{graphics}
\usepackage{amsfonts}
\usepackage{amsmath}

\def\RR{\mathbb{R}}
\def\MM{\mathbb{M}}
\def\OO{\mathbb{O}}
\def\SS{\mathbb{S}}

\def\eps{\varepsilon}
\def\wto{\rightharpoonup}
\def\Aff{{\rm Aff}}
\def\det{{\rm det}}
\def\Z{\mathcal{Z}}

\def\rank{{\rm rank}}

\newtheorem{theorem}{Theorem}[section]
\newtheorem{lemma}[theorem]{Lemma}
\newtheorem{proposition}[theorem]{Proposition}
\newtheorem{corollary}[theorem]{Corollary}

\theoremstyle{definition}
\newtheorem{definition}[theorem]{Definition}

\theoremstyle{remark}
\newtheorem{remark}[theorem]{Remark}

\title{Relaxation theorems in nonlinear elasticity}

\author{Omar Anza Hafsa}
\address{Institut f\"ur Mathematik, Universit\"at Z\"urich, Winterthurerstrasse 190, CH-8057 Z\"urich, Switzerland.}
\email{anza@math.unizh.ch}
\author{Jean-Philippe Mandallena}
\address{``Equipe AVA (Analyse Variationnelle et Applications)" Centre Universitaire de Formation et de Recherche de N\^\i mes, Site des Carmes, Place Gabriel P\'eri - Cedex 01 - 30021 N\^\i mes, France.\newline \indent I3M (Institut de Math\'ematiques et Mod\'elisation de Montpellier) UMR - CNRS 5149, Universit\'e Montpellier II, Place Eug\`ene Bataillon, 34090 Montpellier, France.}
\email{jean-philippe.mandallena@unimes.fr}

\begin{document}

\begin{abstract}
Relaxation theorems which apply to one, two and three-dimensio-nal nonlinear elasticity are proved. We take into account the fact an infinite amount of energy is required to compress a finite line, surface or volume into zero line, surface or volume. However, we do not prevent orientation reversal.
\end{abstract}

\maketitle


\section{Main results}

\subsection{Introduction} Consider an elastic material occupying in a reference configuration $\Omega\subset\RR^N$ with $N=1,2{\rm\ or\ }3$, where $\Omega$ is bounded and open with Lipschitz boundary $\partial\Omega$. The mechanical properties of the material are characterized by a stored-energy function $W:\MM^{3\times N}\to[0,+\infty]$ (assumed to be Borel measurable) in terms of which the total stored-energy is the integral
\begin{equation}\label{Integral}
I(u):=\int_\Omega W(\nabla u(x))dx
\end{equation}
with $\nabla u(x)\in\MM^{3\times N}$ the gradient of $u$ at $x$, where $\MM^{3\times N}$ denotes the space of all real $3\times N$ matrices. In order to take into account the fact that an infinite amount of energy is required to compress a finite line ($N=1$), surface ($N=2$) or volume ($N=3$) into zero line, surface or volume, i.e.,
\begin{equation}\label{DetHypo}
W(\xi)\to+\infty\hbox{ as }\left\{
\begin{array}{ll}|\xi|\to 0&\hbox{if }N=1\\
|\xi_1\land \xi_2|\to 0&\hbox{if }N=2\\
|\det \xi|\to 0&\hbox{if }N=3,\end{array}
\right.
\end{equation}
we consider the following conditions:
\begin{itemize}
\item[(C$_1$)] \em there exist $\alpha,\beta>0$ such that for every $\xi\in\MM^{3\times 1}$,
$$
\hbox{if }|\xi|\geq\alpha\hbox{ then }W(\xi)\leq \beta(1+|\xi|^p)
$$
\end{itemize}
for $N=1$;
\begin{itemize}
\item[(C$_2$)] \em there exist $\alpha,\beta>0$ such that for every $\xi=(\xi_1\mid\xi_2)\in\MM^{3\times 2}$, 
$$
\hbox{if }|\xi_1\land\xi_2|\geq\alpha\hbox{ then }W(\xi)\leq \beta(1+|\xi|^p)
$$
\end{itemize}
for $N=2$, where $\xi_1\land\xi_2$ denotes the cross product of vectors $\xi_1,\xi_2\in\RR^3$; 
\begin{itemize}
\item[(C$_3$)] \em for every $\delta>0$, there exists $c_\delta>0$ such that for every $\xi\in\MM^{3\times 3}$,
$$
\hbox{if }|\det\xi|\geq\delta\hbox{ then }W(\xi)\leq c_\delta(1+|\xi|^p),
$$
\end{itemize}
where $\det\xi$ denotes the determinant of $\xi$, and
\begin{itemize}
\item[(C$_4$)] \em $W(P\xi Q)=W(\xi)$ for all $\xi\in\MM^{3\times 3}$ and all $P,Q\in\SS\OO(3)$
\end{itemize}
for $N=3$, with $\SS\OO(3):=\{Q\in\MM^{3\times 3}:Q^{\rm T}Q=QQ^{\rm T}=I_3\hbox{ and }\det Q=1\}$, where $I_3$ denotes the identity matrix in $\MM^{3\times 3}$ and $Q^{\rm T}$ is the transposed matrix of $Q$. (In fact, (C$_4$) is an additional condition which is not related to (\ref{DetHypo}). However, it means that $W$ is frame-indifferent, i.e., $W(P\xi)=W(\xi)$ for all $\xi\in\MM^{3\times 3}$ and all $P\in\SS\OO(3)$, 
and isotropic, i.e., $W(\xi Q)=W(\xi)$ for all $\xi\in\MM^{3\times 3}$ and all $Q\in\SS\OO(3)$, see for example \cite{hughesmarsden} for more details.)

Fix $p\in]1,+\infty[$, set $W^{1,p}_g(\Omega;\RR^3):=\{u\in W^{1,p}(\Omega;\RR^3):u=g\hbox{ on }\partial\Omega\}$, where  $g$ is given continuous piecewise affine function from $\Omega$ to $\RR^3$, define the integral
$$
\mathcal{Q}I(u):=\int_\Omega \mathcal{Q}W(\nabla u(x))dx,
$$
where $\mathcal{Q}W:\MM^{3\times N}\to[0,+\infty]$ denotes the quasiconvex envelope of $W$, and consider the following assertions: 
\begin{itemize}
\item[(R$_1$)] $\inf\left\{I(u):u\in W^{1,p}_g(\Omega;\RR^3)\right\}=\inf\left\{\mathcal{Q}I(u):u\in W^{1,p}_g(\Omega;\RR^3)\right\}${\rm;}
\item[(R$_2$)] \em if $u_n\wto \overline{u}$ with $\{u_n\}_{n\geq 1}$ minimizing sequence  for $I$ in $W^{1,p}_g(\Omega;\RR^3)$, then $\overline{u}$ is a minimizer for $\mathcal{Q}I$ in $W^{1,p}_g(\Omega;\RR^3)${\rm;}
\item[(R$_3$)] if $\overline{u}$ is a minimizer for $\mathcal{Q}I$ in $W^{1,p}_g(\Omega;\RR^3)$, then there exists a minimizing sequence $\{u_n\}_{n\geq 1}$  for $I$ in $W^{1,p}_g(\Omega;\RR^3)$ such that $u_n\wto \overline{u}$,
\end{itemize}
where ``$\wto$" denotes the weak convergence in $W^{1,p}(\Omega;\RR^3)$. In this paper we prove (see Sect. 1.3) the following relaxation theorems:
\begin{theorem}[$N=1$]\label{1Dtheorem}
If {\rm(C$_1$)} holds and if $W$ is coercive, i.e., $W(\xi)\geq C|\xi|^p$ for all $\xi\in\MM^{3\times N}$ and some $C>0$,  then {\rm(R$_1$)}, {\rm(R$_2$)} and {\rm(R$_3$)} hold.
\end{theorem}
\begin{theorem}[$N=2$]\label{2Dtheorem}
If {\rm(C$_2$)} holds and if $W$ is coercive,  then {\rm(R$_1$)}, {\rm(R$_2$)} and {\rm(R$_3$)} hold.
\end{theorem}
\begin{theorem}[$N=3$]\label{3Dtheorem}
If {\rm(C$_3$)} and {\rm(C$_4$)} hold and if $W$ is coercive,  then {\rm(R$_1$)}, {\rm(R$_2$)} and {\rm(R$_3$)} hold.
\end{theorem}
 Typically, these theorems can be applied with stored-energy functions $W$of the form
$$
W(\xi):=|\xi|^p+\left\{
\begin{array}{ll}
h(|\xi|)&\hbox{if }N=1\\
h(|\xi_1\land\xi_2|)&\hbox{if }N=2\\
h(|\det\xi|)&\hbox{if }N=3,
\end{array}
\right.
$$
for all $\xi\in\MM^{3\times N}$, where $h:[0,+\infty[\to[0,+\infty]$ is Borel measurable and such that for every $\delta>0$, there exists $r_\delta>0$ such that $h(t)\leq r_\delta$ for all $t\geq\delta$ (for example,  $h(0)=+\infty$ and $h(t)=1/t^s$ if $t>0$ with $s>0$).  

\subsection{Outline of the paper} Let $\mathcal{I}:W^{1,p}(\Omega;\RR^3)\to[0,+\infty]$ be defined by
$$
\mathcal{I}(u):=\left\{\begin{array}{cl}
\displaystyle\int_\Omega W(\nabla u(x))dx&\hbox{if }u\in W^{1,p}_g(\Omega;\RR^3)\\
+\infty&\hbox{otherwise,}
\end{array}
\right.
$$
let $\mathcal{Q}\mathcal{I}:W^{1,p}(\Omega;\RR^3)\to[0,+\infty]$ be defined by
$$
\mathcal{Q}\mathcal{I}(u):=\left\{\begin{array}{cl}
\displaystyle\int_\Omega \mathcal{Q}W(\nabla u(x))dx&\hbox{if }u\in W^{1,p}_g(\Omega;\RR^3) \\
+\infty&\hbox{otherwise,}
\end{array}
\right.
$$
and let $\overline{\mathcal{I}}:W^{1,p}(\Omega;\RR^3)\to[0,+\infty]$ be the lower semicontinuous envelope (or relaxed functional) of $\mathcal{I}$ with respect to the weak topology of $W^{1,p}(\Omega;\RR^3)$, i.e.,
$$
\overline{\mathcal{I}}(u):=\inf\left\{\liminf_{n\to+\infty}\mathcal{I}(u_n):u_n\wto u\right\}.
$$
Set $Y:=]0,1[^N$ and $\Aff_0(Y;\RR^3):=\{\phi\in\Aff(Y;\RR^3):\phi=0\hbox{ on }\partial Y\}$, where $\Aff(Y;\RR^3)$ denotes the space of all continuous piecewise affine functions from $Y$ to $\RR^3$, and consider $\Z W:\MM^{3\times N}\to[0,+\infty]$ given by
$$
\Z W(\xi):=\inf\left\{\int_Y W(\xi+\nabla\phi(y))dy:\phi\in \Aff_0(Y;\RR^3)\right\}.
$$
Here is the central theorem of the paper:
\begin{theorem}\label{MainResult} 
If $\Z W$ is of $p$-polynomial growth, i.e., 
$
\Z W(\xi)\leq c(1+|\xi|^p)
$
for all $\xi\in\MM^{3\times N}$ and some $c>0$, then $\overline{\mathcal{I}}=\mathcal{Q}\mathcal{I}$.
\end{theorem}
 Here $m=3$ and $N=1,2\hbox{ or }3$, but the proof of Theorem \ref{MainResult} (given in Section 3) does not depend on the integers $m$ and $N$. This immediately gives the following relaxation result:
\begin{corollary}\label{RelaxMainResult}
Under the hypotheses of Theorem {\rm\ref{MainResult}}, if $W$ is coercive, then {\rm(R$_1$)}, {\rm(R$_2$)} and {\rm(R$_3$)} hold.
\end{corollary}
Such results was proved by Dacorogna in \cite{dacorogna0} when $W$ is continuous and of $p$-polynomial growth. The distinguishing feature here is that Theorem \ref{MainResult} (and so Corollary \ref{RelaxMainResult}) is compatible with (\ref{DetHypo}). More precisely, in Section 4 we prove the following propositions:
\begin{proposition}[$N=1$]\label{1Dth}
If {\rm(C$_1$)} holds then $\Z W$ is of $p$-polynomial growth. 
\end{proposition}
\begin{proposition}[$N=2$]\label{2Dth}
If {\rm(C$_2$)} holds then $\Z W$ is of $p$-polynomial growth.
\end{proposition}
\begin{proposition}[$N=3$]\label{3Dth}
If {\rm(C$_3$)} and {\rm(C$_4$)} hold then $\Z W$ is of $p$-polynomial growth.
\end{proposition}
Theorem \ref{MainResult} follows from Propositions \ref{QuasicFormula} and \ref{Sob=Aff} below whose proofs are given in Section 3:
\begin{proposition}\label{QuasicFormula}
If $\Z W$ is finite then $\mathcal{Q}W=\mathcal{Q}[\Z W]=\Z W$. Furthermore, for $N=1$ we have $\Z W=W^{**}$, where $W^{**}$ denotes the lower semicontinuous convex envelope of $W$.
\end{proposition}
\begin{proposition}\label{Sob=Aff}
$\mathcal{J}_0=\mathcal{J}_1$ with  $\mathcal{J}_0,\mathcal{J}_1:W^{1,p}(\Omega;\RR^3)\to[0,+\infty]$ respectively defined by
$$
\mathcal{J}_0(u):=\inf\left\{\liminf_{n\to+\infty}I(u_n):\Aff_g(\Omega;\RR^3)\ni u_n\wto u\right\}
$$
and
$$
\mathcal{J}_1(u):=\inf\left\{\liminf_{n\to+\infty}\Z I(u_n):\Aff_g(\Omega;\RR^3)\ni u_n\wto u\right\},
$$
where $\Aff_g(\Omega;\RR^3):=\{u\in\Aff(\Omega;\RR^3):u=g\hbox{ on }\partial\Omega\}$ and 
$$
\Z I(u):=\int_\Omega \Z W(\nabla u(x))dx.
$$
\end{proposition}
Taking Proposition \ref{QuasicFormula} into account, from Propositions \ref{1Dth}, \ref{2Dth} and \ref{3Dth}, we see that stored-energy functions $W$ satisfying (C$_1$) for $N=1$, (C$_2$) for $N=2$ and (C$_3$) and (C$_4$) for $N=3$, are not quasiconvex, so that the integral $I(u)$ in (\ref{Integral}) is not weakly lower semicontinuous on $W^{1,p}(\Omega;\RR^3)$ (see \cite[Corollary 3.2]{ballmurat}). Thus, the Direct Method of the Calculus of Variations cannot be applied to study the existence of minimizers of $I$ in $W^{1,p}_g(\Omega;\RR^3)$. For this reason, in the present paper we establish relaxation theorems instead of existence theorems. (In fact, the term ``relaxation" means ``generalized existence",  see \cite{ekeland,daco,buttazzo} for a deeper discussion.)

Other related results can be found in \cite{carde,benbelgacem} where we refer the reader. The present work improves our previous one \cite{anzman2} (see also \cite{anzman4,anzman5}). The main new contribution of the present paper is the treatment of the case $N=3$.

The plan of the paper is as follows. The proofs of Theorems \ref{1Dtheorem}, \ref{2Dtheorem} and \ref{3Dtheorem} are given in Sect. 1.3  (although these can be easily deduced from the previous discussion). Section 2 presents some preliminaries. In Section 3 we prove Propositions \ref{QuasicFormula} and \ref{Sob=Aff} and Theorem \ref{MainResult}. Finally, Section 4 contains the proofs of Propositions \ref{1Dth}, \ref{2Dth} and \ref{3Dth}.

\subsection{Proof of Theorems \ref{1Dtheorem}, \ref{2Dtheorem} and \ref{3Dtheorem}} According to Corollary \ref{RelaxMainResult}, we see that Theorems \ref{1Dtheorem}, \ref{2Dtheorem} and \ref{3Dtheorem} are immediate consequences of respectively Propositions \ref{1Dth}, \ref{2Dth} and \ref{3Dth}.\hfill$\square$


\section{Preliminaries}

In this section we recall some (classical) definitions and results. These will be used throughout the paper.

Let $m,N\geq 1$ be two integers. For any bounded open set $D\subset\RR^N$, we denote by $\Aff(D;\RR^m)$ the space of all continuous piecewise affine functions from $D$ to $\RR^m$, i.e., $u\in\Aff(D;\RR^m)$ if and only if $u$ is continuous and there exists a finite family $(D_i)_{i\in I}$ of open disjoint subsets of $D$ such that $|D\setminus \cup_{i\in I} D_i|=0$ and for every $i\in I$, $\nabla u(x)=\xi_i$ in $D_i$ with $\xi_i\in\MM^{m\times N}$ (where $|\cdot|$ denotes the Lebesgue measure in $\RR^N$).  For any $g\in W^{1,p}(D;\RR^m)$ with $p>1$, we set 
$$
\Aff_g(D;\RR^m):=\{u\in\Aff(D;\RR^m):u=g\hbox{ on }\partial D\}
$$
($\Aff_0(D;\RR^m)$ corresponds to $\Aff_g(D;\RR^m)$ with $g=0$) and 
$$
W^{1,p}_g(D;\RR^m):=\{u\in W^{1,p}(D;\RR^m):u=g\hbox{ on }\partial D\},
$$
(where $\partial D$ denotes the boundary of $D$). Fix $g\in W^{1,p}(\Omega;\RR^m)$ where $\Omega\subset\RR^N$ is bounded and open with Lipschitz boundary. The following density theorem will play an essential role in the proof of Theorem \ref{MainResult}:
 \begin{theorem}[Ekeland-Temam \cite{ekeland}]\label{DensityTheorem}
$\Aff_g(\Omega;\RR^m)$ is dense in $W^{1,p}_g(\Omega;\RR^m)$ with respect to the strong topology of $W^{1,p}(\Omega;\RR^m)$.
 \end{theorem}
 Let $f:\MM^{m\times N}\to[0,+\infty]$ be Borel measurable and let $\Z f:\MM^{m\times N}\to[0,+\infty]$ be defined by
$$
\Z f(\xi):=\inf\left\{\int_Y f(\xi+\nabla\phi(x))dx:\phi\in\Aff_0(Y;\RR^m)\right\}
$$
with $Y:=]0,1[^N$. To prove Propositions \ref{1Dth}, \ref{2Dth} and \ref{3Dth}, we will use  the following properties of $\Z f$:
\begin{proposition}[Fonseca \cite{fonseca}]\label{FonsecaLemma}
\begin{itemize} 
\item[(i)] For every bounded open set $D\subset\RR^N$ with $|\partial D|=0$ and every $\xi\in\MM^{m\times N}$,
$$
\Z f(\xi)=\inf\left\{{1\over |D|}\int_D f(\xi+\nabla\phi(x))dx:\phi\in\Aff_0(D;\RR^m)\right\}.
$$
\item[(ii)] If $\Z f$ is finite then $\Z f$ is rank-one convex, i.e., for every $\xi,\xi^\prime\in\MM^{m\times N}$ with $\rank(\xi-\xi^\prime)\leq 1$, $\Z f(\lambda\xi+(1-\lambda)\xi^\prime)\leq\lambda\Z f(\xi)+(1-\lambda)\Z f(\xi^\prime)$.
 \item[(iii)] If $\Z f$ is finite then $\Z f$ is continuous.
\item[(iv)] For every bounded open set $D\subset\RR^N$ with $|\partial D|=0$, every $\xi\in\MM^{m\times N}$ and every $\phi\in\Aff_0(D;\RR^m)$,
$$
\Z f(\xi)\leq{1\over |D|}\int_D \Z f(\xi+\nabla\phi(x))dx.
$$
\end{itemize}
\end{proposition}
Note that Proposition \ref{FonsecaLemma} is not exactly the one that can found in Fonseca. Nevertheless, it can be proved using the same arguments than the one given by Fonseca (for more details see \cite[Remark A.2]{anzman4}).

Quasiconvexity is the correct concept to deal with multiple integral problems in the Calculus of Variations. For the convenience of the reader, we recall its definition:
\begin{definition}[Morrey \cite{morrey}]
We say that $f$ is quasiconvex if for every $\xi\in\MM^{m\times N}$, every bounded open set $D\subset\RR^N$ with $|\partial D|=0$ and every $\phi\in W^{1,\infty}_0(D;\RR^m)$,
$$
f(\xi)\leq{1\over|D|}\int_D f(\xi+\nabla\phi(x))dx.
$$ 
\begin{remark}\label{ZW=WifQW=W}
Clearly, if $f$ is quasiconvex then $\Z f=f$.
\end{remark}
\end{definition}
By the quasiconvex envelope of $f$, that we denote $\mathcal{Q}f$, we mean the greatest quasiconvex function which less than or equal to $f$. (Thus, $f$ is quasiconvex if and only if $\mathcal{Q}f=f$.) To prove Proposition \ref{QuasicFormula} we will need Theorem \ref{DRTqc}:
\begin{theorem}[Dacorogna \cite{dacorogna0,daco}]\label{DRTqc}
If $f$ is continuous and finite then $\mathcal{Q}f=\Z f$.
\end{theorem}
Let $\mathcal{F}:W^{1,p}(\Omega;\RR^m)\to[0,+\infty]$ be defined by
$$
\mathcal{F}(u):=\left\{\begin{array}{cl}
\displaystyle\int_\Omega f(\nabla u(x))dx&\hbox{if }u\in W^{1,p}_g(\Omega;\RR^m)\\
+\infty&\hbox{otherwise,}
\end{array}
\right.
$$
let $\mathcal{Q}\mathcal{F}:W^{1,p}(\Omega;\RR^m)\to[0,+\infty]$ be given by
$$
\mathcal{Q}\mathcal{F}(u):=\left\{\begin{array}{cl}
\displaystyle\int_\Omega \mathcal{Q}f(\nabla u(x))dx&\hbox{if }u\in W^{1,p}_g(\Omega;\RR^m) \\
+\infty&\hbox{otherwise,}
\end{array}
\right.
$$
and let $\overline{\mathcal{F}}:W^{1,p}(\Omega;\RR^m)\to[0,+\infty]$ be the lower semicontinuous envelope (or relaxed functional) of $\mathcal{I}$ with respect to the weak topology of $W^{1,p}(\Omega;\RR^m)$, i.e.,
$$
\overline{\mathcal{F}}(u):=\inf\left\{\liminf_{n\to+\infty}\mathcal{F}(u_n):u_n\wto u\right\},
$$
where ``$\wto$" denotes the weak convergence in $W^{1,p}(\Omega;\RR^m)$. We close this section with the following integral representation theorem that we will use in the proof of Theorem \ref{MainResult}:
\begin{theorem}[Dacorogna \cite{dacorogna0,daco}]\label{DIRT}
If $f$ is continuous and of $p$-polynomial growth, i.e., $f(\xi)\leq c(1+|\xi|^p)$ for all $\xi\in\MM^{m\times N}$ and some $c>0$, then $\overline{\mathcal{F}}=\mathcal{Q}\mathcal{F}$. 
\end{theorem}


\section{Proof of Propositions \ref{QuasicFormula} and \ref{Sob=Aff} and Theorem \ref{MainResult}}

\subsection{Proof of Proposition \ref{QuasicFormula}} Since $\Z W$ is finite, $\Z W$ is continuous by Proposition \ref{FonsecaLemma}(iii). From Theorem \ref{DRTqc}, we deduce that $\mathcal{Q}[\Z W]=\Z[\Z W]$. But $\Z[\Z W]=\Z W$ by Proposition \ref{FonsecaLemma}(iv), hence $\mathcal{Q}[\Z W]=\Z W$. On the other hand, as $\Z W\leq W$ we have $\mathcal{Q}[\Z W]\leq\mathcal{Q}W$. Moreover, as $\mathcal{Q}W$ is quasiconvex, from Remark \ref{ZW=WifQW=W} we see that $\Z[\mathcal{Q}W]=\mathcal{Q}W$, and consequently $\mathcal{Q}W\leq\mathcal{Q}[\Z W]$. It follows that $\mathcal{Q}W=\mathcal{Q}[\Z W]=\Z W$.

Assume that $N=1$. Then, quasiconvexity is equivalent to convexity (see \cite[Theorem 1.1 p. 102]{daco}). Thus, $\Z W$ is convex (resp. $\Z W^{**}=W^{**}$) since $\Z W=\mathcal{Q}W$ (resp. $W^{**}$ is convex). But $\Z W$ is continuous (resp. $W^{**}\leq W$), hence $\Z W\leq W^{**}$ (resp. $W^{**}\leq\Z W$). It follows that $\Z W=W^{**}$.\hfill$\square$

\subsection{Proof of Proposition \ref{Sob=Aff}} Clearly $\mathcal{J}_1\leq \mathcal{J}_0$. We are thus reduced to prove that 
\begin{equation}\label{J_0=J_1}
\mathcal{J}_0\leq\mathcal{J}_1.
\end{equation}

We need the following lemma:
\begin{lemma}\label{LemmaStep1}
If $u\in\Aff_{g}(\Omega;\RR^3)$ then $\mathcal{J}_0(u)\leq \Z I(u)$.
\end{lemma}
\noindent{\em Proof of Lemma {\rm\ref{LemmaStep1}}.} Let $u\in\Aff_{g}(\Omega;\RR^3)$. By definition, there exists a finite family $(\Omega_i)_{i\in I}$ of open disjoint subsets of $\Omega$ such that $|\Omega\setminus\cup_{i\in I}\Omega_i|=0$ and, for every $i\in I$, $\nabla u(x)=\xi_i$ in $\Omega_i$ with $\xi_i\in\MM^{3\times N}$. Given any $\delta>0$ and any $i\in I$, we consider $\phi_i\in \Aff_0(Y;\RR^3)$ such that
\begin{equation}\label{Zk}
\int_Y W(\xi_i+\nabla\phi_i(y))dy\leq\Z W(\xi_i)+{\delta\over|\Omega|}.
\end{equation}
Fix any integer $n\geq 1$. By Vitali's covering theorem, there exists a finite or countable family $(a_{i,j}+\eps_{i,j}Y)_{j\in J_{i}}$ of disjoint subsets of $\Omega_i$, where $a_{i,j}\in\RR^N$ and $0<\eps_{i,j}<{1\over n}$, such that
$
|\Omega_i\setminus\cup_{j\in J_{i}}(a_{i,j}+\eps_{i,j}Y)|=0
$
 (and so $\sum_{j\in J_i}\eps_{i,j}^N=|\Omega_i|$). Define $\psi_n:\Omega\to\RR^3$ by
$$
\psi_n(x):=
\eps_{i,j}\phi_{i}\left({x-a_{i,j}\over \eps_{i,j}}\right)\hbox{ if }x\in a_{i,j}+\eps_{i,j}Y.
$$
Clearly, for every $n\geq 1$, $\psi_n\in\Aff_0(\Omega;\RR^3)$, $\|\psi_n\|_{L^\infty(\Omega;\RR^3)}\leq {1\over n}\max_{i\in I}\|\phi_i\|_{L^\infty(Y;\RR^3)}$ and $\|\nabla\psi_n\|_{L^\infty(\Omega;\MM^{3\times N})}\leq\max_{i\in I}\|\nabla\phi_i\|_{L^\infty(Y;\MM^{3\times N})}$, hence (up to a subsequence) $\psi_n\stackrel{*}{\wto}0$ in $W^{1,\infty}(\Omega;\RR^3)$, where ``$\stackrel{*}{\wto}$" denotes the weak$^*$ convergence in $W^{1,\infty}(\Omega;\RR^3)$. Consequently,  $\psi_n\wto 0$ in $W^{1,p}(\Omega;\RR^3)$. Moreover, 
\begin{eqnarray*}
I(u+\psi_n)&=&\sum_{i\in I}\int_{\Omega_i}W(\xi_i+\nabla\psi_n(x))dx\\
&=&\sum_{i\in I}\sum_{j\in J_i}\eps_{i,j}^N\int_{Y}W(\xi_i+\nabla\phi_{i}(y))dy\\
&=&\sum_{i\in I}|\Omega_i|\int_{Y}W(\xi_i+\nabla\phi_{i}(y))dy.
\end{eqnarray*}
As $u+\psi_n\in\Aff_{g}(\Omega;\RR^3)$ for all $n\geq 1$ and $u+\psi_n\wto u$ in $W^{1,p}(\Omega;\RR^3)$, from (\ref{Zk}) we deduce that
$$
\mathcal{J}_0(u)\leq\liminf_{n\to+\infty}I(u+\psi_n)\leq\sum_{i\in I}|\Omega_i|\Z W(\xi_i)+\delta=\Z I(u)+\delta,
$$
and the lemma follows.\hfill$\square$

Fix any $u\in W^{1,p}(\Omega;\RR^3)$ and any sequence $u_n\wto u$ with $u_n\in\Aff_{g}(\Omega;\RR^3)$. Using Lemma \ref{LemmaStep1} we have
$
\mathcal{J}_0(u_n)\leq\Z I(u_n)
$
for all $n\geq 1$. Thus,
$$
\mathcal{J}_0(u)\leq\liminf_{n\to+\infty}\mathcal{J}_0(u_n)\leq\liminf_{n\to+\infty}\Z I(u_n),
$$
and (\ref{J_0=J_1}) follows. \hfill $\square$


\subsection{Proof of Theorem \ref{MainResult}} Since $\Z W$ is of $p$-polynomial growth, $\Z W$ is finite, and so $\Z W$ is continuous by Proposition \ref{FonsecaLemma}(iii). From Theorem \ref{DensityTheorem} it follows that
$$
\mathcal{J}_1(u)=\inf\left\{\liminf_{n\to+\infty}\Z I(u_n):W^{1,p}_{g}(\Omega;\RR^3)\ni u_n\wto u\right\}.
$$
But $\mathcal{Q}W=\mathcal{Q}[\Z W]$ by Proposition \ref{QuasicFormula}, hence $\mathcal{J}_1=\mathcal{Q}\mathcal{I}$ by Theorem \ref{DIRT}. On the other hand, given any $u\in W^{1,p}(\Omega;\RR^3)$ and any $u_n\wto u$ with $u_n\in W^{1,p}_{g}(\Omega;\RR^3)$, we have
$
\Z I(u_n)\leq I(u_n)
$
for all $n\geq 1$. Thus,
$$
\mathcal{J}_1(u)\leq\liminf_{n\to+\infty}\Z I(u_n)\leq\liminf_{n\to+\infty}I(u_n),
$$
and so $\mathcal{J}_1\leq\overline{\mathcal{I}}$. But $\overline{\mathcal{I}}\leq\mathcal{J}_0$ and $\mathcal{J}_0=\mathcal{J}_1$ by Proposition \ref{Sob=Aff}, hence $\overline{\mathcal{I}}=\mathcal{J}_1$, and the theorem follows. \hfill$\square$ 


\section{Proof of Propositions \ref{1Dth}, \ref{2Dth} and \ref{3Dth}}


\subsection{Case \boldmath$N=1$\unboldmath} In this section we prove Proposition \ref{1Dth}. 

\medskip

\noindent {\em Proof of Proposition {\rm\ref{1Dth}}.} By (C$_1$) it is clear that if $|\xi|\geq\alpha$ then $\Z W(\xi)\leq\beta(1+|\xi|^p)$. Fix any $\xi\in\MM^{3\times 1}$ such that $|\xi|\leq\alpha$ and consider $\phi\in\Aff_0(Y;\RR^3)$ given by
$$
\phi(x):=\left\{
\begin{array}{ll}
x\nu&\hbox{if }x\in]0,{1\over 2}]\\
(1-x)\nu&\hbox{if }x\in]{1\over 2},1[
\end{array}
\right.
$$
with $\nu\in\RR^3$ such that $|\nu|=2\alpha$. Then,
$$
|\xi+\nabla\phi(x)|=\left\{
\begin{array}{ll}
|\xi+\nu|&\hbox{if }x\in]0,{1\over 2}[\\
|\xi-\nu|&\hbox{if }x\in]{1\over 2},1[,
\end{array}
\right.
$$
hence $|\xi+\nabla\phi(x)|\geq\min\{|\xi+\nu|,|\xi-\nu|\}\geq|\nu|-|\xi|\geq\alpha$ for all $x\in]0,{1\over 2}[\cup]{1\over 2},1[$, and from (C$_1$) we deduce that 
\begin{eqnarray*}
\Z W(\xi)\leq\int_YW(\xi+\nabla\phi(x))dx&\leq&\beta\int_Y(1+|\xi+\nabla\phi(x)|^p)dx\\&\leq&\beta2^{2p}\max\{1,\alpha^p\}(1+|\xi|^p).
\end{eqnarray*} 
It follows that $\Z W(\xi)\leq\beta2^{2p}\max\{1,\alpha^p\}(1+|\xi|^p)$ for all $\xi\in\MM^{3\times 1}$. \hfill$\square$


\subsection{Case \boldmath$N=2$\unboldmath} In this section we prove Proposition \ref{2Dth}. We begin with the following lemma.
\begin{lemma}\label{Lemma1forN=2}
Under {\rm(C$_2$)} there exists $\gamma>0$ such that for all $\xi=(\xi_1\mid\xi_2)\in\MM^{3\times 2}$, 
$$
\hbox{if $\min\{|\xi_1+\xi_2|,|\xi_1-\xi_2|\}\geq\alpha$ then $\Z W(\xi)\leq\gamma(1+|\xi|^p).$}
$$
\end{lemma}
\begin{proof}
Let $\xi=(\xi_1\mid\xi_2)\in\MM^{3\times 2}$ be such that $\min\{|\xi_1+\xi_2|,|\xi_1-\xi_2|\}\geq\alpha$ (with $\alpha>0$ given by (C$_2$)). Then, one the three possibilities holds:
\begin{itemize}
\item[(i)] $|\xi_1\land\xi_2|\not=0$;
\item[(ii)] $|\xi_1\land\xi_2|=0$ with $\xi_1\not=0$;
\item[(iii)] $|\xi_1\land\xi_2|=0$ with $\xi_2\not=0$.
\end{itemize}
Set 
$
D:=\{(x_1,x_2)\in\RR^2:x_1-1<x_2<x_1+1\hbox{ and }-x_1-1<x_2<1-x_1\}
$ 
and, for each $t\in\RR$, define $\varphi_t\in\Aff_0(D;\RR)$ by
$$
\varphi_t(x_1,x_2):=\left\{
\begin{array}{ll}
-tx_1+t(x_2+1)&\hbox{if }(x_1,x_2)\in\Delta_1\\
t(1-x_1)-tx_2&\hbox{if }(x_1,x_2)\in\Delta_2\\
tx_1+t(1-x_2)&\hbox{if }(x_1,x_2)\in\Delta_3\\
t(x_1+1)+tx_2&\hbox{if }(x_1,x_2)\in\Delta_4
\end{array}
\right.
$$
with
\begin{itemize} 
\item[]$
\Delta_1:=\{(x_1,x_2)\in D:x_1\geq 0\hbox{ and } x_2\leq 0\};
$
\item[]$
\Delta_2:=\{(x_1,x_2)\in D:x_1\geq 0\hbox{ and }x_2\geq 0\};
$
\item[]$
\Delta_3:=\{(x_1,x_2)\in D:x_1\leq 0\hbox{ and }x_2\geq 0\};
$
\item[]$
\Delta_4:=\{(x_1,x_2)\in D:x_1\leq 0\hbox{ and }x_2\leq 0\}.
$
\end{itemize}
Consider $\phi\in\Aff_0(D;\RR^3)$ given by 
$$
\phi:=(\varphi_{\nu_1},\varphi_{\nu_2},\varphi_{\nu_3})\hbox{ with }\left\{\begin{array}{ll}\nu={\xi_1\land \xi_2\over|\xi_1\land \xi_2|}&\hbox{if (i) is satisfied}\\
|\nu|=1 \hbox{ and }\langle\xi_1,\nu\rangle=0&\hbox{if (ii) is satisfied}\\
|\nu|=1 \hbox{ and }\langle\xi_2,\nu\rangle=0&\hbox{if (iii) is satisfied}\end{array}\right.
$$
($\nu_1$, $\nu_2$, $\nu_3$ are the components of the vector $\nu$). Then, 
$$
\xi+\nabla\phi(x)=\left\{
\begin{array}{ll}
(\xi_1-\nu\mid \xi_2+\nu)&\hbox{if }x\in{\rm int}(\Delta_1)\\
(\xi_1-\nu\mid \xi_2-\nu)&\hbox{if }x\in{\rm int}(\Delta_2)\\
(\xi_1+\nu\mid \xi_2-\nu)&\hbox{if }x\in{\rm int}(\Delta_3)\\
(\xi_1+\nu\mid \xi_2+\nu)&\hbox{if }x\in{\rm int}(\Delta_4)
\end{array}
\right.
$$
(where ${\rm int}(E)$ denotes the interior of the set $E$). Taking Proposition \ref{FonsecaLemma}(i) into account, it follows that 
\begin{eqnarray}\label{Z_1}
\Z W(\xi)&\leq&{1\over 4}\Big(W(\xi_1-\nu\mid \xi_2+\nu)+W(\xi_1-\nu\mid \xi_2-\nu)\\
&&+\ W(\xi_1+\nu\mid \xi_2-\nu)+W(\xi_1+\nu\mid \xi_2+\nu)\Big).\nonumber
\end{eqnarray}
But
$
|(\xi_1-\nu)\land(\xi_2+\nu)|^2=|\xi_1\land \xi_2+(\xi_1+\xi_2)\land\nu|^2
=|\xi_1\land \xi_2|^2+|(\xi_1+\xi_2)\land\nu|^2
\geq|(\xi_1+\xi_2)\land\nu|^2,
$
and so 
$$
|(\xi_1+\nu)\land(\xi_2-\nu)|\geq |(\xi_1+\xi_2)\land\nu|=|\xi_1+\xi_2|.
$$
Similarly, we obtain: 
\begin{itemize}
\item[] $|(\xi_1-\nu)\land(\xi_2-\nu)|\geq |\xi_1-\xi_2|$;
\item[] $|(\xi_1+\nu)\land(\xi_2-\nu)|\geq |\xi_1+\xi_2|$;
\item[] $|(\xi_1+\nu)\land(\xi_2+\nu)|\geq |\xi_1-\xi_2|$.
\end{itemize}
Thus, $|(\xi_1-\nu)\land(\xi_2+\nu)|\geq\alpha$, $|(\xi_1-\nu)\land(\xi_2-\nu)|\geq\alpha$, $|(\xi_1+\nu)\land(\xi_2-\nu)|\geq\alpha$ and $|(\xi_1+\nu)\land(\xi_2+\nu)|\geq\alpha$, because $\min\{|\xi_1+\xi_2|,|\xi_1-\xi_2|\}\geq\alpha$. Using (C$_2$) it follows that 
\begin{eqnarray*}
W(\xi_1-\nu\mid \xi_2+\nu)&\leq&\beta(1+|(\xi_1-\nu\mid \xi_2+\nu)|^p)\\
&\leq&\beta 2^p(1+|(\xi_1\mid \xi_2)|^p+|(-\nu\mid \nu)|^p)\\
&\leq&\beta 2^{2p+1}(1+|\xi|^p).
\end{eqnarray*}
In the same manner, we have: 
\begin{itemize}
\item[] $W(\xi_1-\nu\mid \xi_2-\nu)\leq \beta 2^{2p+1}(1+|\xi|^p)$;
\item[] $W(\xi_1+\nu\mid \xi_2-\nu)\leq \beta 2^{2p+1}(1+|\xi|^p)$;
\item[] $W(\xi_1+\nu\mid \xi_2+\nu)\leq \beta 2^{2p+1}(1+|\xi|^p)$,
\end{itemize}
and, from (\ref{Z_1}), we conclude that
$
\Z W(\xi)\leq \beta 2^{2p+1}(1+|\xi|^p).
$
\end{proof}

\medskip 

\noindent{\em Proof of Proposition {\rm\ref{2Dth}}.} Let $\xi=(\xi_1\mid \xi_2)\in\MM^{3\times 2}$. Then, one the four possibilities holds:
\begin{itemize}
\item[(i)] $|\xi_1\land\xi_2|\not=0$;
\item[(ii)] $|\xi_1\land\xi_2|=0$ with $\xi_1=\xi_2=0$;
\item[(iii)] $|\xi_1\land\xi_2|=0$ with $\xi_1\not=0$;
\item[(iv)] $|\xi_1\land\xi_2|=0$ with $\xi_2\not=0$.
\end{itemize}
For each $t\in\RR$, define $\varphi_t\in\Aff_0(Y;\RR)$ by 
$$
\varphi_t(x_1,x_2):=\left\{
\begin{array}{ll}
tx_2&\hbox{if }(x_1,x_2)\in\Delta_1\\
t(1-x_1)&\hbox{if }(x_1,x_2)\in\Delta_2\\
t(1-x_2)&\hbox{if }(x_1,x_2)\in\Delta_3\\
tx_1&\hbox{if }(x_1,x_2)\in\Delta_4
\end{array}
\right.
$$
with
\begin{itemize} 
\item[]$
\Delta_1:=\big\{(x_1,x_2)\in Y:x_2\leq x_1\leq -x_2+1\big\};
$
\item[]$
\Delta_2:=\big\{(x_1,x_2)\in Y:-x_1+1\leq x_2\leq x_1\big\};
$
\item[]$
\Delta_3:=\big\{(x_1,x_2)\in Y:-x_2+1\leq x_1\leq x_2\big\};
$
\item[]$
\Delta_4:=\big\{(x_1,x_2)\in Y:x_1\leq x_2\leq -x_1+1\big\}.
$
\end{itemize}
Consider $\phi\in\Aff_0(Y;\RR^3)$ given by 
$$
\phi:=\big(\varphi_{\nu_1},\varphi_{\nu_2},\varphi_{\nu_3}\big)\hbox{ with }
\left\{
\begin{array}{ll}
\nu={\alpha(\xi_1\land \xi_2)\over|\xi_1\land \xi_2|}&\hbox{if (i) is satisfied}\\  
|\nu|=\alpha&\hbox{if (ii) is satisfied}\\
|\nu|=\alpha\hbox{ and }\langle\xi_1,\nu\rangle=0& \hbox{if (iii) is satisfied}\\
|\nu|=\alpha\hbox{ and }\langle\xi_2,\nu\rangle=0& \hbox{if (iv) is satisfied}
\end{array}
\right.
$$
($\nu_1$, $\nu_2$, $\nu_3$ are the components of the vector $\nu$ and $\alpha>0$ is given by (C$_2$)). Then, 
$$
\xi+\nabla\phi(x)=\left\{
\begin{array}{ll}
(\xi_1\mid \xi_2+\nu)&\hbox{if }x\in{\rm int}(\Delta_1)\\
(\xi_1-\nu\mid \xi_2)&\hbox{if }x\in{\rm int}(\Delta_2)\\
(\xi_1\mid \xi_2-\nu)&\hbox{if }x\in{\rm int}(\Delta_3)\\
(\xi_1+\nu\mid \xi_2)&\hbox{if }x\in{\rm int}(\Delta_4)
\end{array}
\right.
$$
(where ${\rm int}(E)$ denotes the interior of the set $E$). Taking Proposition \ref{FonsecaLemma}(iv) into account, it follows that
\begin{eqnarray}\label{ZzZ}
\Z W(\xi)&\leq&{1\over 4}\Big(\Z W(\xi_1\mid \xi_2+\nu)+\Z W(\xi_1-\nu\mid \xi_2)\\
&&+\ \Z W(\xi_1\mid \xi_2-\nu)+\Z W(\xi_1+\nu\mid \xi_2)\Big).\nonumber
\end{eqnarray}
But
$
|\xi_1+(\xi_2+\nu)|^2=|(\xi_1+\xi_2)+\nu|^2
=|\xi_1+\xi_2|^2+|\nu|^2
=|\xi_1+\xi_2|^2+\alpha^2
\geq \alpha^2,
$
hence
$
|\xi_1+(\xi_2+\nu)|\geq \alpha.
$
Similarly, we obtain
$
|\xi_1-(\xi_2+\nu)|\geq \alpha,
$
and so
$$
\min\{|\xi_1+(\xi_2+\nu)|,|\xi_1-(\xi_2+\nu)|\}\geq \alpha.
$$
In the same manner, we have: 
\begin{itemize}
\item[] $\min\{|(\xi_1-\nu)+\xi_2|,|(\xi_1-\nu)-\xi_2|\}\geq\alpha$;
\item[] $\min\{|\xi_1+(\xi_2-\nu)|,|\xi_1-(\xi_2-\nu)|\}\geq\alpha$;
\item[] $\min\{|(\xi_1+\nu)+\xi_2|,|(\xi_1+\nu)-\xi_2|\}\geq\alpha$. 
\end{itemize}
Using Lemma \ref{Lemma1forN=2} it follows that
\begin{eqnarray*}
\Z W(\xi_1\mid \xi_2+\nu)&\leq&\gamma(1+|(\xi_1\mid \xi_2+\nu)|^p)\\
&\leq&\gamma 2^p(1+|(\xi_1\mid\xi_2)|^p+|(0\mid\nu)|^p)\\
&\leq&\max\{1,\alpha^p\}\gamma2^{p+1}(1+|\xi|^p).
\end{eqnarray*}
Similarly, we obtain:
\begin{itemize}
\item[] $\Z W(\xi_1-\nu\mid \xi_2)\leq\max\{1,\alpha^p\}\gamma2^{p+1}(1+|\xi|^p)$;
\item[] $\Z W(\xi_1\mid \xi_2-\nu)\leq\max\{1,\alpha^p\}\gamma2^{p+1}(1+|\xi|^p)$;
\item[] $\Z W(\xi_1+\nu\mid \xi_2)\leq\max\{1,\alpha^p\}\gamma2^{p+1}(1+|\xi|^p)$,
\end{itemize} 
and, from (\ref{ZzZ}), we conclude  that 
$
\Z W(\xi)\leq\max\{1,\alpha^p\}\gamma2^{p+1}(1+|\xi|^p).
$ 
\hfill$\square$


\subsection{Case \boldmath$N=3$\unboldmath} In this section we prove Proposition \ref{3Dth}. We begin with three lemmas.
\begin{lemma}\label{Lemma1forN=3}
If {\rm(C$_3$)} holds then $\Z W$ is finite.
\end{lemma}
\begin{proof}
Clearly, if $\xi\in\MM^{3\times 3}_*$ then $\Z W(\xi)<+\infty$ with $\MM^{3\times 3}_*:=\{\xi\in\MM^{3\times 3}:\det\xi\not=0\}$. We are thus reduced to prove that  $\Z W(\xi)<+\infty$ for all $\xi\in\MM^{3\times 3}\setminus\MM^{3\times 3}_*$. Fix $\xi=(\xi_1\mid\xi_2\mid\xi_3)\in\MM^{3\times 3}\setminus\MM^{3\times 3}_*$ where $\xi_1,\xi_2,\xi_3\in\RR^3$ are the columns of $\xi$. Then $\rank(\xi)\in\{0,1,2\}$ (where $\rank(\xi)$ denotes the rank of the matrix $\xi$).

\vskip0.5mm

{\em Step }1. {\em We prove that if $\rank(\xi)=2$ then $\Z W(\xi)<+\infty$}. Without loss of generality we can assume that there exist $\lambda,\mu\in\RR$ such that $\xi_3=\lambda\xi_1+\mu\xi_2$.  Given any $s\in\RR^*$, consider $D\subset\RR^3$ given by
\begin{equation}\label{BoundedOpenSet}
D:={\rm int}(\cup_{i=1}^8\Delta_{i}^s)
\end{equation}
(where ${\rm int}(E)$ denotes the interior of the set $E$) with:
\begin{itemize}
\item[] $\Delta_{1}^s:=\{(x_1,x_2,x_3)\in \RR^3:x_1\geq 0,\;x_2\geq 0,\;x_3\geq 0 \hbox{ and }x_1+x_2+sx_3\leq 1\}$;
\item[] $\Delta_{2}^s:=\{(x_1,x_2,x_3)\in \RR^3:x_1\leq 0,\;x_2\geq 0,\;x_3\geq 0 \hbox{ and }-x_1+x_2+sx_3\leq 1\}$;
\item[] $\Delta_{3}^s:=\{(x_1,x_2,x_3)\in \RR^3:x_1\leq 0,\;x_2\leq 0,\;x_3\geq 0 \hbox{ and }-x_1-x_2+sx_3\leq 1\}$;
\item[] $\Delta_{4}^s:=\{(x_1,x_2,x_3)\in \RR^3:x_1\geq 0,\;x_2\leq 0,\;x_3\geq 0 \hbox{ and }x_1-x_2+sx_3\leq 1\}$;
 \item[] $\Delta_{5}^s:=\{(x_1,x_2,x_3)\in\RR^3:x_1\geq 0,\;x_2\geq 0,\;x_3\leq 0 \hbox{ and }x_1+x_2-sx_3\leq 1\}$;
 \item[] $\Delta_{6}^s:=\{(x_1,x_2,x_3)\in \RR^3:x_1\leq 0,\;x_2\geq 0,\;x_3\leq 0 \hbox{ and }-x_1+x_2-sx_3\leq 1\}$;
 \item[] $\Delta_{7}^s:=\{(x_1,x_2,x_3)\in \RR^3:x_1\leq 0,\;x_2\leq 0,\;x_3\leq 0 \hbox{ and }-x_1-x_2-sx_3\leq 1\}$;
 \item[] $\Delta_{8}^s:=\{(x_1,x_2,x_3)\in \RR^3:x_1\geq 0,\;x_2\leq 0,\;x_3\leq 0 \hbox{ and }x_1-x_2-sx_3\leq 1\}$;
 \end{itemize}
 Clearly, $D$ is bounded, open and $|\partial D|=0$. For each $t\in\RR$, define $\varphi_{s,t}\in\Aff_0(D;\RR)$ by
\begin{equation}\label{AffineFunction}
\varphi_{s,t}(x_1,x_2,x_3):=\left\{
\begin{array}{ll}
-t(x_1+1)-tx_2-tsx_3&\hbox{if }(x_1,x_2,x_3)\in\Delta_{1}^s\\
tx_1-t(x_2+1)-tsx_3&\hbox{if }(x_1,x_2,x_3)\in\Delta_{2}^s\\
t(x_1-1)+tx_2-tsx_3&\hbox{if }(x_1,x_2,x_3)\in\Delta_{3}^s\\
-tx_1+t(x_2-1)-tsx_3&\hbox{if }(x_1,x_2,x_3)\in\Delta_{4}^s\\
-t(x_1+1)-tx_2+tsx_3&\hbox{if }(x_1,x_2,x_3)\in\Delta_{5}^s\\
tx_1-t(x_2+1)+tsx_3&\hbox{if }(x_1,x_2,x_3)\in\Delta_{6}^s\\
t(x_1-1)+tx_2+tsx_3&\hbox{if }(x_1,x_2,x_3)\in\Delta_{7}^s\\
-tx_1+t(x_2-1)+tsx_3&\hbox{if }(x_1,x_2,x_3)\in\Delta_{8}^s.
\end{array}
\right.
\end{equation}
(Note that $\varphi_{s,t}$ is simply the only function, affine on every $\Delta_i^s$, null on $\partial D$, such that $\varphi_{s,t}(0)=-t$.) Fix $s\in\RR^*\setminus\{\lambda-\mu,-(\lambda-\mu),\lambda+\mu,-(\lambda+\mu)\}$ and consider $\phi\in\Aff_0(D;\RR^3)$ given by
$$
\phi:=(\varphi_{s,\nu_1},\varphi_{s,\nu_2},\varphi_{s,\nu_3})\hbox{ with }\nu:={\xi_1\land\xi_2\over|\xi_1\land\xi_2|^2},
$$
($\nu_1$, $\nu_2$, $\nu_3$ are the components of the vector $\nu$). Then,
$$
\xi+\nabla\phi(x)=\left\{
\begin{array}{ll}
(\xi_1-\nu\mid\xi_2-\nu\mid \xi_3-s\nu)&\hbox{if }x\in{\rm int}(\Delta_{1}^s)\\
(\xi_1+\nu\mid \xi_2-\nu\mid \xi_3-s\nu)&\hbox{if }x\in {\rm int}(\Delta_{2}^s)\\
(\xi_1+\nu\mid \xi_2+\nu\mid \xi_3-s\nu)&\hbox{if }x\in {\rm int}(\Delta_{3}^s)\\
(\xi_1-\nu\mid \xi_2+\nu\mid \xi_3-s\nu)&\hbox{if }x\in {\rm int}(\Delta_{4}^s)\\
(\xi_1-\nu\mid\xi_2-\nu\mid \xi_3+s\nu)&\hbox{if }x\in {\rm int}(\Delta_{5}^s)\\
(\xi_1+\nu\mid \xi_2-\nu\mid \xi_3+s\nu)&\hbox{if }x\in {\rm int}(\Delta_{6}^s)\\
(\xi_1+\nu\mid \xi_2+\nu\mid \xi_3+s\nu)&\hbox{if }x\in {\rm int}(\Delta_{7}^s)\\
(\xi_1-\nu\mid \xi_2+\nu\mid \xi_3+s\nu)&\hbox{if }x\in {\rm int}(\Delta_{8}^s).
\end{array}
\right.
$$
As $\det\xi=0$, $\xi_1\land\xi_3=\mu(\xi_1\land\xi_2)$ and $\xi_2\land\xi_3=\lambda(\xi_2\land\xi_1)$ we have
$$
|\det(\xi+\nabla\phi(x))|=\left\{
\begin{array}{ll}
|s+(\lambda-\mu)|&\hbox{if }x\in{\rm int}(\Delta_{1}^s)\cup{\rm int}(\Delta_{7}^s)\\
|s-(\lambda+\mu)|&\hbox{if }x\in{\rm int}(\Delta_{2}^s)\cup{\rm int}(\Delta_{8}^s)\\
|s-(\lambda-\mu)|&\hbox{if }x\in{\rm int}(\Delta_{3}^s)\cup{\rm int}(\Delta_{5}^s)\\
|s+(\lambda+\mu)|&\hbox{if }x\in{\rm int}(\Delta_{4}^s)\cup{\rm int}(\Delta_{6}^s).\\
\end{array}
\right.
$$
It follows that for a.e. $x\in D$,
$
|\det(\xi+\nabla\phi(x)|\geq\min\{|s+(\lambda-\mu)|,|s-(\lambda+\mu)|,|s-(\lambda-\mu)|,|s+(\lambda+\mu)|\}=:\delta
$
($\delta>0$). Taking Proposition \ref{FonsecaLemma}(i) into account and using  {\rm(${\rm C}_3$)}, we see that there exists $c_\delta>0$ such that 
$$
\Z W(\xi)\leq{1\over|D|}\int_DW(\xi+\nabla\phi(x))\leq c_\delta+{c_\delta\over|D|}\|\xi+\nabla \phi\|^p_{L^p(D;\RR^3)},
$$
which implies that $\Z W(\xi)<+\infty$.

\vskip0.5mm

{\em Step }2. {\em We prove that if $\rank(\xi)=1$ then $\Z W(\xi)<+\infty$}. Without loss of generality we can assume that there exist $\lambda,\mu\in\RR$ such that $\xi_2=\lambda\xi_1$ and $\xi_3=\mu\xi_1$. Consider $D\subset\RR^3$ given by (\ref{BoundedOpenSet}) with $s\in\RR^*\setminus\{-\mu,\mu\}$, and define $\phi\in\Aff_0(D;\RR^3)$ by $\phi:=(\varphi_{s,\nu_1},\varphi_{s,\nu_2},\varphi_{s,\nu_3})$ with $\nu=(\nu_1,\nu_2,\nu_3)\in\RR^3\setminus\{0\}$ such that $\langle\nu,\xi_1\rangle=0$, where, for every $i\in\{1,2,3\}$, $\varphi_{s,\nu_i}$ is defined by (\ref{AffineFunction}) with $t=\nu_i$. By Proposition \ref{FonsecaLemma}(iv) we have
\begin{eqnarray*}
\Z W(\xi)&\leq& {1\over 8}\Big(\Z W(\xi_1-\nu\mid\xi_2-\nu\mid \xi_3-s\nu)+\Z W(\xi_1+\nu\mid\xi_2-\nu\mid \xi_3-s\nu)\\
&&\ \ +\Z W(\xi_1+\nu\mid\xi_2+\nu\mid \xi_3-s\nu)+\Z W(\xi_1-\nu\mid\xi_2+\nu\mid \xi_3-s\nu)\\
&&\ \ +\Z W(\xi_1-\nu\mid\xi_2-\nu\mid \xi_3+s\nu)+\Z W(\xi_1+\nu\mid\xi_2-\nu\mid \xi_3+s\nu)\\
&&\ \ +\Z W(\xi_1+\nu\mid\xi_2+\nu\mid \xi_3+s\nu)+\Z W(\xi_1-\nu\mid\xi_2+\nu\mid \xi_3+s\nu)\Big).
\end{eqnarray*}
Noticing that $s\in\RR^*\setminus\{-\mu,\mu\}$ it is easy to see that:
\begin{itemize}
\item[] $\rank(\xi_1-\nu\mid\xi_2-\nu\mid \xi_3-s\nu)=2$;
\item[] $\rank(\xi_1+\nu\mid\xi_2-\nu\mid \xi_3-s\nu)=2$;
\item[] $\rank(\xi_1+\nu\mid\xi_2+\nu\mid \xi_3-s\nu)=2$;
\item[] $\rank(\xi_1-\nu\mid\xi_2+\nu\mid \xi_3-s\nu)=2$;
\item[] $\rank(\xi_1-\nu\mid\xi_2-\nu\mid \xi_3+s\nu)=2$;
\item[] $\rank(\xi_1+\nu\mid\xi_2-\nu\mid \xi_3+s\nu)=2$;
\item[] $\rank(\xi_1+\nu\mid\xi_2+\nu\mid \xi_3+s\nu)=2$;
\item[] $\rank(\xi_1-\nu\mid\xi_2+\nu\mid \xi_3+s\nu)=2$,
\end{itemize}
and using Step  1 we deduce that $\Z W(\xi)<+\infty$.

\vskip0.5mm

{\em Step }3. {\em We prove that $\Z W(0)<+\infty$}. This follows from Step 2 by using Proposition \ref{FonsecaLemma}(iv) with $D\subset\RR^3$ given by (\ref{BoundedOpenSet}) with $s\in\RR^*$, and $\phi\in\Aff_0(D;\RR^3)$ defined by $\phi:=(\varphi_{s,\nu_1},\varphi_{s,\nu_2},\varphi_{s,\nu_3})$ with $(\nu_1,\nu_2,\nu_3)\in\RR^3\setminus\{0\}$, where, for every $i\in\{1,2,3\}$, $\varphi_{s,\nu_i}$ is defined by (\ref{AffineFunction}) with $t=\nu_i$.
\end{proof}
\begin{lemma}\label{Lemma2forN=3}
Under {\rm(C$_3$)} there exists $c>0$ such that for every $\xi\in\MM^{3\times 3}$,
$$
\hbox{if }\xi\hbox{ is diagonal then }\Z W(\xi)\leq c(1+|\xi|^p). 
$$
\end{lemma}
\begin{proof}
Combining Lemma \ref{Lemma1forN=3} with Proposition \ref{FonsecaLemma}(ii), we deduce that $\Z W$ is continuous, and so there exists $c_0>0$ such that $\Z W(\xi)\leq c_0$ for all $\xi\in\MM^{3\times 3}$ with $|\xi|^2\leq 3$. Moreover, it is obvious that $\Z W(\xi)\leq c_1(1+|\xi|^p)$ for all $\xi\in\MM^{3\times 3}$ such that  $|\det\xi|\geq 1$, where $c_1>0$ is given by (C$_3$) with $\delta=1$. We are thus led to consider $\xi\in\MM^{3\times 3}$ such that $\xi$ is diagonal, $|\det\xi|\leq 1$ and $|\xi|^2\geq 3$, i.e., $\xi=(\xi_{ij})$ with $\xi_{ij}=0$ if $i\not=j$, $|\xi_{11}\xi_{22}\xi_{33}|\leq 1$ and $|\xi_{11}|^2+|\xi_{22}|^2+|\xi_{33}|^2\geq 3$. Then, one the six possibilities holds:
\begin{itemize}
\item[(i)] $|\xi_{11}|\leq 1$, $|\xi_{22}|\geq 1$ and $|\xi_{33}|\geq 1$;
\item[(ii)] $|\xi_{22}|\leq 1$, $|\xi_{33}|\geq 1$ and $|\xi_{11}|\geq 1$;
\item[(iii)] $|\xi_{33}|\leq 1$, $|\xi_{11}|\geq 1$ and $|\xi_{22}|\geq 1$;
\item[(iv)] $|\xi_{11}|\geq 1$, $|\xi_{22}|\leq 1$ and $|\xi_{33}|\leq 1$;
\item[(v)] $|\xi_{22}|\geq 1$, $|\xi_{33}|\leq 1$ and $|\xi_{11}|\leq 1$;
\item[(vi)] $|\xi_{33}|\geq 1$, $|\xi_{11}|\leq 1$ and $|\xi_{22}|\leq 1$.
\end{itemize}

{\em Claim }1. {\em There exists $c_2>0$ such that if $\xi$ is diagonal with $|\det\xi|\leq 1$ and satisfies either {\rm(i)}, {\rm(ii)} or {\rm(iii)}, then $\Z W(\xi)\leq c_2(1+|\xi|^p)$.} Consider $D\subset\RR^3$ given by (\ref{BoundedOpenSet}) with $s=1$, and define $\phi\in\Aff_0(D;\RR^3)$ by $\phi:=(\varphi_{1,\nu_1},\varphi_{1,\nu_2},\varphi_{1,\nu_3})$, where 
$$
(\nu_1,\nu_2,\nu_3)=\left\{
\begin{array}{ll}
(2,0,0)&\hbox{if (i) is satisfied}\\
(0,2,0)&\hbox{if (ii) is satisfied}\\
(0,0,2)&\hbox{if (iii) is satisfied,}
\end{array}
\right.
$$
and, for every $i\in\{1,2,3\}$, $\varphi_{1,\nu_i}$ is defined by (\ref{AffineFunction}) with $s=1$ and $t=\nu_i$. It is then easy to see that for a.e. $x\in D$,
$$
|\det(\xi+\nabla\phi(x))|\geq\left\{
\begin{array}{ll}
|2|\xi_{22}||\xi_{33}|-|\det\xi||&\hbox{if (i) is satisfied}\\
|2|\xi_{11}||\xi_{33}|-|\det\xi||&\hbox{if (ii) is satisfied}\\
|2|\xi_{11}||\xi_{22}|-|\det\xi||&\hbox{if (iii) is satisfied,}
\end{array}
\right.
$$
so that  $|\det(\xi+\nabla\phi(x))|\geq 1$. Taking Proposition \ref{FonsecaLemma}(i) into account, using (C$_3$) and noticing that $|\nabla\phi(x)|=2\sqrt{3}$ for a.e. $x\in D$, we deduce that 
$
\Z W(\xi)\leq c_2(1+|\xi|^p)
$
with $c_2:=c_12^p(1+(2\sqrt{3})^p)$.

\vskip0.5mm

{\em Claim }2. {\em There exists $c_3>0$ such that if $\xi$ is diagonal and satisfies either {\rm(iv)}, {\rm(v)} or {\rm(vi)}, then $\Z W(\xi)\leq c_3(1+|\xi|^p)$.} Let $\zeta\in\MM^{3\times 3}$ be a rank-one diagonal matrix defined by
$$
\zeta_{11}:=\left\{\begin{array}{ll}\xi_{22}+{\rm sign(\xi_{22})}&\hbox{if (iv) is satisfied}\\0&\hbox{if either (v) or (vi) is satisfied;}\end{array}\right.
$$
$$
\zeta_{22}:=\left\{\begin{array}{ll}\xi_{33}+{\rm sign(\xi_{33})}&\hbox{if (v) is satisfied}\\0&\hbox{if either (iv) or (vi) is satisfied;}\end{array}\right.
$$
$$
\zeta_{33}:=\left\{\begin{array}{ll}\xi_{11}+{\rm sign(\xi_{11})}&\hbox{if (vi) is satisfied}\\0&\hbox{if either (iv) or (v) is satisfied,}\end{array}\right.
$$
where ${\rm sign}(r)=1$ if $r\geq 0$ and ${\rm sign}(r)=-1$ if $r<0$. Then, $\xi^+:=\xi+\zeta$ and $\xi^{-}:=\xi-\zeta$ are diagonal matrices such that: 
\begin{itemize}
\item[(a)] $|\xi^+_{11}|\geq 1$, $|\xi^+_{22}|\geq 1$, $|\xi^+_{33}|\leq 1$ and $|\xi^-_{11}|\geq 1$, $|\xi^-_{22}|\geq 1$, $|\xi^-_{33}|\leq 1$ if (iv) is satisfied;
\item[(b)] $|\xi^+_{11}|\leq 1$, $|\xi^+_{22}|\geq 1$, $|\xi^+_{33}|\geq 1$ and $|\xi^-_{11}|\leq 1$, $|\xi^-_{22}|\geq 1$, $|\xi^-_{33}|\geq 1$ if (v) is satisfied;
\item[(c)] $|\xi^+_{11}|\geq 1$, $|\xi^+_{22}|\leq 1$, $|\xi^+_{33}|\geq 1$ and $|\xi^-_{11}|\geq 1$, $|\xi^-_{22}|\leq 1$, $|\xi^-_{33}|\geq 1$ if (vi) is satisfied.
\end{itemize}
Combining Lemma \ref{Lemma1forN=3} with Proposition \ref{FonsecaLemma}(ii) we deduce that $\Z W$ is rank-one convex, so that
\begin{equation}\label{Rank-OneConvexity}
\Z W(\xi)\leq {1\over 2}\big(\Z W(\xi^+)+\Z W(\xi^-)\big).
\end{equation}
According to (a), (b) and (c), from Claim 1 we see that  $\Z W(\xi^+)\leq c_2(1+|\xi^+|^p)$ (resp. $\Z W(\xi^-)\leq c_2(1+|\xi^+|^p)$) if $|\det \xi^+|\leq 1$ (resp. $|\det \xi^-|\leq 1$). On the other hand,  by (C$_3$) we have  $\Z W(\xi^+)\leq c_1(1+|\xi^+|^p)$ (resp. $\Z W(\xi^-)\leq c_1(1+|\xi^+|^p)$) if $|\det \xi^+|\geq 1$ (resp. $|\det \xi^-|\geq 1$). Noticing that $|\xi^+|^p\leq 2^{2p}(1+|\xi|^p)$ (resp. $|\xi^-|^p\leq 2^{2p}(1+|\xi|^p)$) and using (\ref{Rank-OneConvexity}), we deduce that $\Z W(\xi)\leq c_3(1+|\xi|^p)$ with $c_3:=2^{2p}\max\{c_1,c_2\}$.

\vskip0.5mm

From Claims 1 and 2, it follows that for every $\xi\in\MM^{3\times 3}$, if $\xi$ is diagonal with $|\xi|^2\geq 3$ and $|\det\xi|\leq 1$ then $\Z W(\xi)\leq c_4(1+|\xi|^p)$ with $c_4:=\max\{c_2,c_3\}$. Setting $c:=\max\{c_0,c_4\}$ we conclude that $\Z W(\xi)\leq c(1+|\xi|^p)$ for all $\xi\in\MM^{3\times 3}$ such that $\xi$ is diagonal.
\end{proof}
\begin{lemma}\label{Lemma3forN=3}
If {\rm(C$_4$)} holds then $\Z W(P\xi Q)=\Z W(\xi)$ for all $\xi\in\MM^{3\times 3}$ and all $P, Q\in\SS\OO(3)$.
\end{lemma}
\begin{proof}
It is suffices to show that 
\begin{itemize}
\item[(i)] $\Z W(P\xi Q)\leq \Z W(\xi)\hbox{ for all }\xi\in\MM^{3\times 3}\hbox{ and all }P,Q\in\SS\OO(3)$.
\end{itemize}
Indeed, given $\xi\in\MM^{3\times 3}$ and $P,Q\in\SS\OO(3)$, we have $\xi=P^{\rm T}(P\xi Q)Q^{\rm T}$, and using (i) we obtain $\Z W(\xi)\leq\Z W(P\xi Q)$. Moreover, (i) is equivalent to 
\begin{itemize}
\item[(ii)] $\Z W(P\xi)\leq \Z W(\xi)\hbox{ for all }\xi\in\MM^{3\times 3}\hbox{ and all }P\in\SS\OO(3)$
\end{itemize}
and
\begin{itemize}
\item[(iii)] $\Z W(\xi Q)\leq \Z W(\xi)\hbox{ for all }\xi\in\MM^{3\times 3}\hbox{ and all }Q\in\SS\OO(3)$.
\end{itemize} 
Indeed, (ii) (resp. (iii)) follows from (i) with $Q=I_3$ (resp. $P=I_3$), where $I_3$ is the identity matrix in $\MM^{3\times 3}$. On the other hand, given $\xi\in\MM^{3\times 3}$ and  $P,Q\in\SS\OO(3)$,  by (ii) (resp. (iii)) we have $\Z W(P(\xi Q))\leq\Z W(\xi Q)$ (resp. $\Z W(\xi Q)\leq\Z W(\xi)$), and so $\Z W(P\xi Q)\leq\Z W(\xi)$. We are thus reduced to prove (ii) and (iii). 

{\em Proof of {\rm(ii)}.} Fix any $\phi\in \Aff_0(Y;\RR^3)$ and set $\varphi:=P\phi$. Then, $\varphi\in \Aff_0(Y;\RR^3)$ and $\nabla\varphi=P\nabla\phi$, hence
\begin{eqnarray*}
\Z W(P\xi)\leq\int_Y W(P\xi+\nabla\varphi(x))dx&=&\int_Y W\big(P(\xi+P^{\rm T}\nabla\varphi(x))\big)dx\\
&=&\int_YW\big(P(\xi+\nabla\phi(x))\big)dx.
\end{eqnarray*}
From (C$_4$) we deduce that
$$
\Z W(P\xi)\leq\int_YW(\xi+\nabla\phi(x))dx
$$
for all $\phi\in\Aff_0(Y;\RR^3)$, which implies that 
$
\Z W(P\xi)\le \Z W(\xi).
$

{\em Proof of {\rm(iii)}.} By Vitali's covering theorem, there exists a finite or countable family $(a_i+\eps_i Q^{\rm T}Y)_{i\in I}$ of disjoint subsets of $Y$, where $a_i\in\RR^3$ and $0<\eps_i<1$, such that $|Y\setminus\cup_{i\in I}(a_i+\eps_i Q^{\rm T}Y)|=0$ (and so $\sum_{i\in I}\eps_i^3=1$). Fix any $\phi\in \Aff_0(Y;\RR^3)$ and define $\varphi\in\Aff_0(Y;\RR^3)$ by
$$
\varphi(x)=\eps_i\phi\left(Q{{x-a_i}\over{\eps_i}}\right)\hbox{ if } x\in a_i+\eps_iQ^{\rm T}Y.
$$
Then,
\begin{eqnarray*}
\Z W(\xi Q)\leq\int_Y W(\xi Q+\nabla\varphi(x))dx&=&\sum_{i\in I}\eps_i^3\int_{Y} W\left(\xi Q+\nabla\phi\left(x\right)Q\right)dx\\
&=&\int_YW\big((\xi+\nabla\phi(x))Q\big)dx.
\end{eqnarray*}
From (C$_4$) we deduce that
$$
\Z W(\xi Q)\leq\int_YW(\xi+\nabla\phi(x))dx
$$
for all $\phi\in\Aff_0(Y;\RR^3)$, which implies that $\Z W(\xi Q)\le \Z W(\xi)$.
\end{proof}

\noindent{\em Proof of Proposition {\rm\ref{3Dth}}.} Fix any $\xi\in\MM^{3\times 3}_*$ (with $\MM^{3\times 3}_*:=\{\xi\in\MM^{3\times 3}:\det\xi\not=0\}$) and consider $P\in\SS\OO(3)$ given by $P:=\xi M^{-1}$ with
$$
M:=\left\{
\begin{array}{ll}
\sqrt{\xi^{\rm T}\xi}&\hbox{if }\det\xi>0\\
-\sqrt{\xi^{\rm T}\xi}&\hbox{if }\det\xi<0.
\end{array}
\right.
$$
As $M$ is symmetric, there exist $Q\in\SS\OO(3)$ and  $\zeta\in\MM^{3\times 3}$ such that $\zeta$ is diagonal and $M=Q^{\rm T}\zeta Q$, hence
$
\xi=PQ^{\rm T}\zeta Q.
$
Consequently, $\Z W(\xi)=\Z W(\zeta)$ by Lemma \ref{Lemma3forN=3}. Noticing that $|\zeta|=|\xi|$, from Lemma \ref{Lemma2forN=3} we deduce that there exists $c>0$ such that $\Z W(\xi)\leq c(1+|\xi|^p)$ for all $\xi\in\MM^{3\times 3}_*$. Combining Lemma \ref{Lemma1forN=3} with Proposition \ref{FonsecaLemma}(iii), we see that $\Z W$ is continuous, and using the fact that $\MM^{3\times 3}_*$ is dense in $\MM^{3\times 3}$, we conclude that  $\Z W(\xi)\leq c(1+|\xi|^p)$ for all $\xi\in\MM^{3\times 3}$. \hfill$\square$



\begin{thebibliography}{}
\bibitem{anzman2} {\sc Anza Hafsa, O., Mandallena, J.-P.}: {\em Relaxation of variational problems in two-dimensional nonlinear elasticity}, to appear on Ann. Mat. Pura Appl. 
\bibitem{anzman4} {\sc Anza Hafsa, O., Mandallena, J.-P.}: {\em The nonlinear membrane energy{\rm:} variational derivation under the constraint {\rm``$\det\nabla u\not=0$"}}, to appear on J. Math. Pures Appl.
\bibitem{anzman5} {\sc Anza Hafsa, O., Mandallena, J.-P.}: {\em The nonlinear membrane energy{\rm:} variational derivation under the constraint {\rm``$\det\nabla u>0$"}}, submitted
\bibitem{ballmurat} {\sc Ball, J. M., Murat, F.}: {\em $W^{1,p}$-quasiconvexity and variational problems for multiple integrals}, J. Funct. Anal. {\bf 58}, 225-253 (1984)
\bibitem{benbelgacem}{\sc Ben Belgacem, H.}: {\em Relaxation of singular functionals defined on Sobolev spaces}, ESAIM Control Optimal Calc. Var.  {\bf 5}, 71-85 (2000)
\bibitem{buttazzo} {\sc Buttazzo, G.}: Semicontinuity, relaxation and integral representation problems in the Calculus of Variations, Pitman Res. Longman, Harlow (1989) 
\bibitem{carde} {\sc Carbone, L., De Arcangelis, R.}: {Unbounded functionals in the Calculus of Variations{\rm:} representation, relaxation and homogenization}, Chapman \& Hall/CRC (2001)
\bibitem{dacorogna0}
{\sc Dacorogna, B.}: {\em Quasiconvexity and relaxation of nonconvex problems in the Calculus of Variations}, J. Funct. Anal. {\bf 46}, 102-118 (1982)
\bibitem{daco} {\sc Dacorogna, B.}: Direct methods in the Calculus of Variations, Springer, Berlin (1989) 
\bibitem{ekeland} {\sc Ekeland, I., Temam, R.}: Analyse convexe et probl\`emes variationnels, Dunod, GauthierVillars, Paris (1974) 
\bibitem{fonseca} {\sc Fonseca, I.}: {\em The lower quasiconvex envelope of the stored energy function for an elastic crystal}, J. Math. Pures et Appl. {\bf 67}, 175-195 (1988)
\bibitem{hughesmarsden} {\sc Marsden, J. E., Hughes, T. J. R.}: Mathematical foundations of elasticity, Prentice-Hall (1983)
\bibitem{morrey} {\sc Morrey, C. B.}: {\em Quasiconvexity and lower semicontinuity of multiple integrals}, Pacific J. Math. {\bf 2}, 25-53 (1952)
\end{thebibliography}
\end{document}